\documentclass[letterpaper, 10 pt, conference]{ieeeconf}  % Comment this line out if you need a4paper

\IEEEoverridecommandlockouts                              % This command is only needed if 
\usepackage{amsmath} % assumes amsmath package installed
\usepackage{amssymb}  % assumes amsmath package installed
\usepackage{mathtools}  % assumes amsmath package installed
\usepackage{amsthm}
\usepackage{algorithm}
\usepackage{algpseudocode}

\newtheorem{theorem}{Theorem}[section]
\newtheorem{proposition}[theorem]{Proposition}

\newtheorem{assumption}[theorem]{Assumption}
\newtheorem{remark}[theorem]{Remark}

\usepackage{xcolor}
\usepackage{comment}
\newif\ifcommentandrea
\commentandreatrue

\usepackage{subcaption}
\usepackage{siunitx}
\usepackage{hyperref}

\title{\LARGE \bf
Inherently robust suboptimal MPC \\for autonomous racing with anytime feasible SQP
}

\author{Logan Numerow, Andrea Zanelli, Andrea Carron and Melanie N. Zeilinger
\thanks{The authors are with the Institute for Dynamical Systems and Control, ETH Zurich, ZH-8092, Switzerland: {\tt\small lnumerow@student.ethz.ch, [zanellia|carrona|mzeilinger]@ethz.ch}}
}

\begin{document}

\maketitle
\thispagestyle{empty}
\pagestyle{empty}

%%%%%%%%%%%%%%%%%%%%%%%%%%%%%%%%%%%%%%%%%%%%%%%%%%%%%%%%%%%%%%%%%%%%%%%%%%%%%%%%
\begin{abstract}
    In recent years, the increasing need for high-performance controllers in applications like autonomous driving has motivated the development of optimization routines tailored to specific control problems. In this paper, we propose an efficient inexact model predictive control (MPC) strategy for autonomous miniature racing with inherent robustness properties. We rely on a feasible sequential quadratic programming (SQP) algorithm capable of generating feasible intermediate iterates such that the solver can be stopped after any number of iterations, without jeopardizing recursive feasibility. In this way, we provide a strategy that computes suboptimal and yet feasible solutions with a computational footprint that is much lower than state-of-the-art methods based on the computation of locally optimal solutions.
    Under suitable assumptions on the terminal set and on the controllability properties of the system, we can state that, for any sufficiently small disturbance affecting the system's dynamics, recursive feasibility can be guaranteed.
    We validate the effectiveness of the proposed strategy in simulation and by deploying it onto a physical experiment with autonomous miniature race cars. Both the simulation and experimental results demonstrate that, using the feasible SQP method, a feasible solution can be obtained with moderate additional computational effort compared to strategies that resort to early termination without providing a feasible solution. At the same time, the proposed method is significantly faster than the state-of-the-art solver Ipopt.
\end{abstract}

%%%%%%%%%%%%%%%%%%%%%%%%%%%%%%%%%%%%%%%%%%%%%%%%%%%%%%%%%%%%%%%%%%%%%%%%%%%%%%%%
\section{Introduction}

Model predictive control (MPC) is an optimization-based control strategy that is especially suited to address difficult control problems including, e.g., nonlinearities and constraints, or large scale problems \cite{Rawlings2017}. Due to the computational burden associated with the solution of the underlying nonconvex programs, MPC has been originally applied, starting from the late 70s, to systems with slow dynamics such as in the process and chemical industry. However, due to the considerable algorithmic progress made in the last decades and - potentially even more prominently - due to the drastic increase in computational power of embedded control units, it has progressively become a viable approach for systems with faster dynamics. Among other fields, MPC is finding its way into autonomous driving applications, where it is gradually becoming an established strategy for trajectory planning and control. For example, the works in \cite{Verschueren2014a}, \cite{Liniger2014}, \cite{Kabzan2020}, \cite{Kloeser2020} report on successful application of MPC for autonomous driving to small and large scale cars, while companies are starting to commercialize planners for autonomous driving based on MPC, e.g. \cite{embotech}.
Despite the growing industrial impact of MPC, one of the main obstacles to its successful deployment remains the computation time needed by underlying solvers.
\par
In the context of safety critical applications such as autonomous racing, it is natural to prioritize constraint satisfaction over optimality. However, most existing solutions do not offer the possibility to decouple (local) optimality from feasibility. Most state-of-the-art approaches to numerically solve the underlying nonconvex programs can be categorized into two groups: \textit{a)} computationally expensive \textit{"exact"} methods that calculate locally optimal, and therefore feasible, solutions, and \textit{b)} typically faster \textit{"inexact"} methods that provide approximate solutions, often deemed infeasible. For instance, the methods used in \cite{Verschueren2014a}, \cite{Liniger2014}, and \cite{Kloeser2020} belong to the second group, as they rely on early stopping of a sequential quadratic programming (SQP) strategy without any guarantees on the feasibility of the approximate solution.

\begin{figure}
    \centering
    \vspace{0.3cm}
    \includegraphics[width=0.485\textwidth]{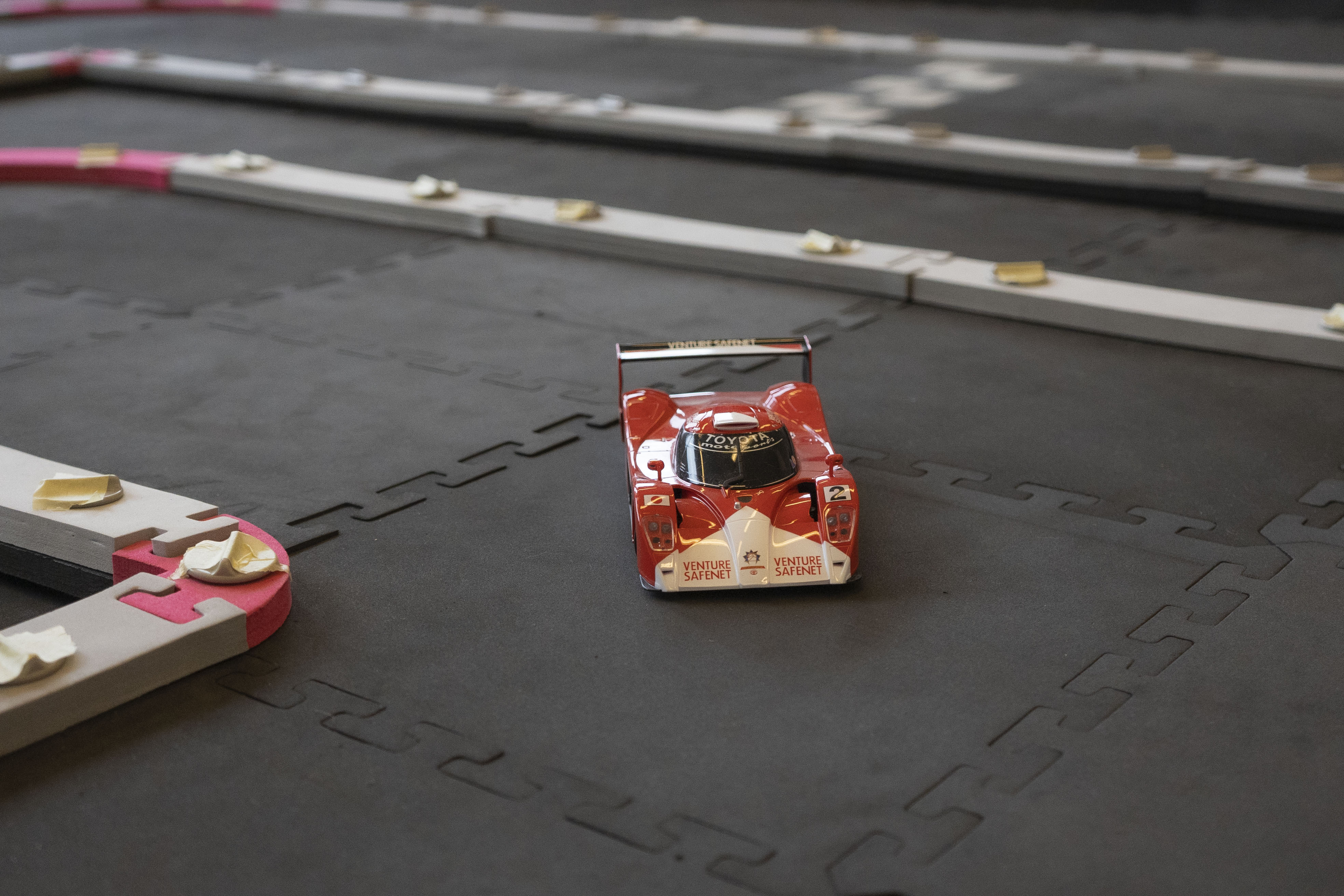}
    \caption{Kyosho Mini-Z 1:28 scale miniature race car.}
    \label{fig:crs_car}
\end{figure}

\subsection{Contribution}
In order to reduce the computational footprint of MPC for autonomous racing while retaining feasibility of the computed trajectories, we propose a strategy based on a feasible SQP algorithm from \cite{zanelli2021}. The algorithm has the advantage that its intermediate iterates are feasible solutions; the algorithm can be terminated after any number of iterations without incurring constraint violation. The contributions of the paper are twofold:
\begin{enumerate}
    \item We demonstrate how the algorithm can be applied to MPC for autonomous racing. We propose a control strategy that, under mild assumptions, is inherently robust to disturbances affecting the system dynamics, and we achieve recursive feasibility by designing a time-varying terminal invariant set based on the offline computation of a feasible periodic trajectory over a full~lap.
    \item We validate the proposed method through both simulation and hardware experiments using miniature race cars. In both scenarios, the anytime feasible SQP algorithm reliably computes feasible trajectories in only slightly increased computation time compared to state-of-the-art approaches that, however, do not provide feasible solutions and can therefore cause safety-critical constraint violations.
\end{enumerate}
The paper is structured as follows. In Section \ref{sec:systemmodel} and~\ref{sec:ocp}, we introduce the dynamical model of the miniature RC car and we describe the optimal control formulation that defines the control policy, respectively. Section \ref{sec:fsqp} describes the feasible SQP method. In Section \ref{sec:rec_feas}, we present the main theoretical result of the paper which shows $M$-step open loop recursive feasibility for the proposed approach, an extension of recursive feasibility which is inherently robust to solver failures and delays. Finally, Section \ref{sec:experiments} describes the simulation and experimental results.
\subsection{Notation}
Throughout the paper, we will denote the Euclidean norm  by $\| \cdot \|$, when referring to vectors,
and, with the same notation, to the spectral norm
$\| A \| \vcentcolon = \sqrt{\lambda_{\max}\left(A^{\top}A\right)}$
when referring
to a real matrix $A$. All
vectors are column vectors and we denote the concatenation of two vectors
by
\begin{equation}
    (x,y)\vcentcolon=\begin{bmatrix}x \\ y\end{bmatrix}.
\end{equation}
We denote the gradient of any function by
$\nabla f(x) = \frac{\partial f}{\partial x}(x)^{\top}$ and the Euclidean ball
of radius $r$ centered at $x$ as $\mathcal{B}(x,r) \vcentcolon=\{ y \,\vcentcolon \,\| x - y \| \leq r\}$.
For some positive integer $n$, we denote by $\mathbb{S}^n_{++}$ ($\mathbb{S}^n_+$) the set of symmetric positive definite (semidefinite) matrices of dimension $n$.
For a dynamical system $x_+ = \psi(x,u,w)$ with state $x \in \mathbb{R}^{n_x}$, input $u \in \mathbb{R}^{n_u}$ and disturbance $w \in \mathbb{R}^{n_w}$, we denote the $k$-times composition of $\psi$, i.e., the result of a $k$-steps forward time simulation, by $\phi(x,u,w;k) \vcentcolon \mathbb{R}^{n_x} \times \mathbb{R}^{n_u \cdot (k-1)} \times \mathbb{R}^{n_w \cdot (k-1)} \rightarrow \mathbb{R}^{n_x}$.
Finally, we denote the identity matrix by $\mathbb{I}$.

\section{System Model}\label{sec:systemmodel}

In this work, we consider a miniature RC car modeled using a dynamic bicycle model formulation \cite{Rajamani2006, Liniger2014}. The state of the model is $x = (p_x, p_y, \gamma, v_f, v_l, \omega, \tau, \delta, \theta)$, where~$(p_x, p_y)$ is the position and $\gamma$ the heading angle of the car in the global coordinate frame, $(v_f, v_l)$ are forward and lateral (left-positive) velocities in the body frame, $\omega$ is the angular velocity, $\tau$ and $\delta$ are the torque and steering angle, and $\theta$ the progress along the track measured by arclength along the center line. The control input is $u = (d\tau, d\delta, d\theta)$, i.e. the time derivatives of $\tau$, $\delta$ and $\theta$. Note that the true system inputs~$\tau$ and $\delta$ are treated as model states, and their derivatives as control inputs. This is done in order to be able to penalize and constrain changes in the system inputs, which can significantly reduce jittering and improve smoothness of the closed-loop trajectories.

The evolution of the system is described by the differential equations

\begin{align}
    \label{eq:eom}
    \begin{split}
        \dot{x} =
        \begin{bmatrix}
            v_f \cos(\gamma) - v_l \sin(\gamma)         \\
            v_f \sin(\gamma) + v_l \cos(\gamma)         \\
            \omega                                      \\
            m^{-1}(F_x - F_f \sin(\delta)) + v_l \omega \\
            m^{-1}(F_r + F_f \cos(\delta)) - v_f \omega \\
            I_z^{-1}(F_f l_f \cos(\delta) - F_r l_r)    \\
            d\tau                                       \\
            d\delta                                     \\
            d\theta                                     \\
        \end{bmatrix},
    \end{split}
\end{align}

where $m$ is the car mass, $I_z$ is the yaw moment of inertia, and $l_{f/r}$ are the distances between the center of mass and the front and rear axles, respectively. The lateral tire forces~$F_{f/r}$ are modeled using a Pacejka tire model \cite{Pacejka}
\begin{equation}
    \begin{aligned}
         &  & \alpha_f & = -\arctan\left(v_l v_f^{-1}\right) - l_f\omega v_f^{-1} + \delta, \\
         &  & \alpha_r & = -\arctan\left(v_l v_f^{-1}\right) + l_r\omega v_f^{-1},          \\
         &  & F_{f/r}  & = D_{f/r}\sin(C_{f/r}\arctan(B_{f/r}\alpha_{f/r})),
    \end{aligned}
\end{equation}
where $\alpha_{f/r}$ are the tire slip angles, and the longitudinal force is modeled as a polynomial function of torque input and velocity:
$F_x = (C_{m1} - C_{m2} v_f) \tau - C_{d} v_f^2 - C_{roll}.$

The exact discrete-time model is denoted by $x_{i+1} = \psi(x_i, u_i, w_i)$, where $w_i$ is a disturbance. This is approximated using fourth-order Runge-Kutta (RK4) integration in the MPC controller and in simulation. In the MPC design we use a nominal model with $w_i = 0, i=0,\dots,N-1$.

\section{Optimal control problem formulation}\label{sec:ocp}
%\subsection{System Constraints}

In this section, we describe the underlying optimal control problem defining the MPC feedback policy.
\begin{comment}
The optimization problem is given in \eqref{eq:nmpc}. Constraints describe (in order) state evolution according to the dynamics model, limits on state and input values, soft contouring constraints, an initial condition at the current measured state and a terminal set constraint (see section \ref{sub:terminalset} for details on terminal set).
\end{comment}
\begin{comment}
\begin{alignat}{4}
    \nonumber
    \min_{\substack{\mathbf{u}_0, \mathbf{u}_1, \cdots, \mathbf{u}_{N-1}                                           \\\mathbf{y}_0, \mathbf{y}_1, \cdots, \mathbf{y}_N}} \quad & f(\mathbf{x}) && \\
    \nonumber
    \text{s.t.} \quad & \mathbf{y}_{i+1} = \mathbf{y}_i + g(\mathbf{y}_i, \mathbf{u}_i), \quad &  & i=0\dots (N-1) \\
    \label{eq:mpc}
                      & \mathbf{y}_i \in \mathcal{Y}, \quad                                    &  & i=0\dots N     \\
    \nonumber
                      & \mathbf{u}_i \in \mathcal{U}, \quad                                    &  & i=0\dots (N-1) \\
    \nonumber
                      & |\langle x_i-\bar{x}_i, y_i-\bar{y}_i \rangle| < r, \quad              &  & i=0\dots N     \\
    \nonumber
                      & \mathbf{y}_0 = \mathbf{y}(k)                                           &  &
\end{alignat}

This is a trajectory optimization problem with $N$ shooting nodes at which optimal control inputs are computed. In controlling the vehicle, only the first control input $\mathbf{u}_0$ is applied before reoptimizing with a shifted horizon.
\end{comment}
%\subsection{State Constraints}
The system is subject to linear input constraints
\begin{equation}
    \begin{aligned}
        d\tau_{min}   & \leq & d\tau   & \leq & d\tau_{max},   \\
        d\delta_{min} & \leq & d\delta & \leq & d\delta_{max}, \\
    \end{aligned}
\end{equation}
denoted $(d\tau, d\delta) \in \mathcal{U}$ and linear state constraints only on the physical inputs of the system
\begin{equation}
    \begin{aligned}
        \tau_{min}   & \leq & \tau   & \leq & \tau_{max},   \\
        \delta_{min} & \leq & \delta & \leq & \delta_{max}, \\
    \end{aligned}
\end{equation}
denoted $(\tau, \delta) \in \mathcal{X}$. Physical inputs are not affected by disturbances, so these can be satisfied robustly as hard state constraints.

Nonlinear state constraints which ensure that the vehicle remains within track boundaries are given by
\begin{equation}
    \label{eq:stateconstraints}
    \pi(x) = \left(p_x-\bar{p}_x(\theta)\right)^2 + \left(p_y-\bar{p}_y(\theta)\right)^2 - \left(\frac{D}{2}\right)^2 \leq \xi,
\end{equation}
where point $(\bar{p}_x(\theta), \bar{p}_y(\theta))$ is the position along the track center line at arclength $\theta$ measured from the start position,~$D$ is the track width and $\xi$ is a slack variable. The position of the vehicle is therefore (softly) constrained to a circle around the corresponding point on the center line with diameter equal to the track width. %(see Figure \ref{fig:contouringconstraints}).  
State constraints are implemented using a cubic spline approximation of $\bar{p}_x(\theta)$ and $\bar{p}_y(\theta)$, with control points every 3 \unit{\centi\metre} along the track.

% \begin{figure}
%     \centering
%     \includegraphics[width=0.485\textwidth]{Figures/contouring_constraints_v1.png}
%     \caption{An illustration of the state constraints used by the MPC controller. From top to bottom: a trajectory satisfying the constraints, a trajectory violating the constraints by leaving the track, and a trajectory violating the constraints without leaving the track.}
%     \label{fig:contouringconstraints}
% \end{figure}

%\subsection{Optimization Objective}

The stage cost is defined as \eqref{eq:stagecost}
\begin{equation}
    \label{eq:stagecost}
    \begin{aligned}
         &  & l(x, u) = & (q_Ce_C(x))^2 + (q_Le_L(x))^2 + (q_{d\tau}d\tau)^2             \\
         &  &           & + (q_{d\delta}d\delta)^2 + (q_{d\theta}(d\theta - \bar{v}))^2,
    \end{aligned}
\end{equation}
with cost weights $q_C$, $q_L$, $q_{d\tau}$, $q_{d\delta}$ and $q_{d\theta}$, and target speed $\bar{v}$, adapted from the model predictive contouring control cost function from \cite{Liniger2014}.

Contouring cost $e_C$ and lag cost $e_L$ correspond to cross-track (lateral) and along-track (longitudinal) deviation from $(\bar{p}_x(\theta)$, $\bar{p}_y(\theta))$ respectively. These are computed using a linearization of the track center line about the initial guess $\hat{\theta}_i$ of the progress state variable, i.e.,

\begin{align}
    \label{eq:statecost}
    \begin{split}
        & e_{C,i} = \frac{d\bar{p}_y}{d\theta}\bigg|_{\hat{\theta}_i} \left(p_{x,i} - \bar{p}_{x,i}(\hat{\theta}_i) - \frac{d\bar{p}_x}{d\theta}\bigg|_{\hat{\theta}i} (\theta_i - \hat{\theta}_i)\right) \\
        & \hspace{4ex} - \frac{d\bar{p}_x}{d\theta}\bigg|_{\hat{\theta}_i} \left(p_{y,i} - \bar{p}_{y,i}(\hat{\theta}_i) - \frac{d\bar{p}_y}{d\theta}\bigg|_{\hat{\theta}_i} (\theta_i - \hat{\theta}_i)\right), \\
        & e_{L,i} = \frac{d\bar{p}_x}{d\theta}\bigg|_{\hat{\theta}_i} \left(p_{x,i} - \bar{p}_{x,i}(\hat{\theta}_i) - \frac{d\bar{p}_x}{d\theta}\bigg|_{\hat{\theta}_i} (\theta_i - \hat{\theta}_i)\right) \\
        & \hspace{4ex} + \frac{d\bar{p}_y}{d\theta}\bigg|_{\hat{\theta}_i} \left(p_{y,i} - \bar{p}_{y,i}(\hat{\theta}_i) - \frac{d\bar{p}_y}{d\theta}\bigg|_{\hat{\theta}_i} (\theta_i - \hat{\theta}_i)\right). \\
    \end{split}
\end{align}

The resulting optimal control problem \eqref{eq:nmpc} reads

\begin{subequations}\label{eq:nmpc}
    \begin{align}
         & \underset{\begin{subarray}{c}
                             x_0, \dots, x_N \\
                             u_0, \dots, u_{N-1} \\
                             \xi_0, \dots, \xi_{N-1}
                         \end{subarray}}{\min} &  & \sum_{i=0}^{N-1} \left(l(x_i, u_i) + \mu \cdot \xi_i\right)               \\
         & \,\,\,\quad \rm{s.t.}         &  & x_0 - \tilde{x} = 0                                                         \\
         &                               &  & \hspace{-0.7ex}\left.\begin{aligned} & x_{i+1} - \psi(x_i,u_i,0) = 0        \\
                & (\tau_i, \delta_i) \in \mathcal{X}   \\
                & (d\tau_i, d\delta_i) \in \mathcal{U} \\
                & \pi(x_i) \leq \xi_i                  \\
                & \xi_i \geq 0                         \\
                & i = 0,\ldots , N-1                   \\
                                                                   \end{aligned}\right\} & \\
         &                               &  & \pi_{N,t}(x_N) \leq 0. \label{eq:nmpc_terminal_constraint}
    \end{align}
\end{subequations}

If the parameter $\mu>0$ is sufficiently large, we can guarantee that problem \eqref{eq:nmpc} has the same solution as its counterpart with hard constraints \cite{Nocedal2006}. Moreover, the terminal constraint~(\ref{eq:nmpc_terminal_constraint}) corresponds to a parameterization of the time-varying terminal set $\mathcal{X}_f(t)$, given by
\begin{equation}
    \mathcal{X}_f(t) \vcentcolon = \{x \,\vcentcolon \pi_{N,t}(x_N) \leq 0\},
\end{equation}
which will be described in detail in Section \ref{sec:rec_feas}. In particular,
this setting with soft state constraints and a hard terminal constraint is similar to, e.g., the one in \cite{Allan2017}, and allows us to establish inherent robust recursive feasibility.
We will refer to problem \eqref{eq:nmpc} as $\mathcal{P}(\tilde{x},t)$ since it depends on the parameter~$\tilde{x}$ representing the current state of the system and on the time-varying terminal constraint defined by $\pi_{N,t}$.

\section{Feasible sequential quadratic programming}\label{sec:fsqp}
In the following, we describe the numerical method used to solve \eqref{eq:nmpc}.
Introducing
\begin{equation}
    y\vcentcolon=\left( x_0, \dots x_N, u_0, \dots, u_{N-1}, \xi_0, \dots, \xi_{N-1}\right)
\end{equation} to refer
to the primal variables, we can rewrite \eqref{eq:nmpc} in the compact form
\begin{equation}\label{eq:compactp}
    \begin{aligned}
         & \underset{y}{\min} &  & f(y)         \\
         & \,\,\rm{s.t.}      &  & g(y) = 0     \\
         &                    &  & h(y) \leq 0,
    \end{aligned}
\end{equation}
where $y \in \mathbb{R}^n$, $g \vcentcolon \mathbb{R}^{n} \rightarrow
    \mathbb{R}^{n_g}$ and $h \vcentcolon \mathbb{R}^{n} \rightarrow
    \mathbb{R}^{n_h}$.
\par
With this compact notation, the QP subproblem associated with a standard sequential quadratic programming (SQP) problem reads as
\begin{equation}\label{eq:compact_sqp}
    \begin{aligned}
         & \underset{\begin{subarray}{c}
                             \Delta y
                         \end{subarray}}{\min}     \!\!\! \!\!
         &                                       & \nabla
        f(\hat{y})^{\top}\Delta y + \frac{1}{2}\Delta y^{\top}M\Delta y         \\
         & \,\,\, \rm{s.t.}                      &        & g(\hat{y}) + \nabla
        g(\hat{y})^{\top} \Delta y = 0,                                         \\
         &                                       &        & h(\hat{y}) + \nabla
        h(\hat{y})^{\top} \Delta y \leq 0.                                      \\
    \end{aligned}
\end{equation}
Each SQP iteration consists of solving \ref{eq:compact_sqp} and using the result to update the solution estimate $\hat y$. Here, we have introduced~$M$ to denote the chosen approximation of the Hessian of the Lagrangian, and $\lambda$ and $\mu$ are the Lagrange multipliers associated with
the equality and inequality constraints, respectively. Following \cite{zanelli2021}, we instead solve modified subproblems, where we do not recompute $\nabla f(\hat{y})$, $\nabla g(\hat{y})$, $\nabla h(\hat{y}), M$ for each QP solve, but use instead fixed derivatives evaluated at the current so-called \textit{outer} iteration~$\tilde{y}$. That is, at each \textit{inner} iteration, we solve the following QP:
\begin{equation}\label{eq:compact_zo_sqp}
    \begin{aligned}
         & \underset{\Delta y}{\min} &  & a(\tilde{y}, \hat{y})^{\top} \Delta y + \frac{1}{2}\Delta y^{\top} M \Delta y \\
         & \,\,\rm{s.t.}             &  & g(\hat{y})  + \nabla g(\tilde{y})^{\top} \Delta y = 0,                        \\
         &                           &  & h(\hat{y}) + \nabla h(\tilde{y})^{\top} \Delta y \leq 0,
    \end{aligned}
\end{equation}
where $a(\tilde{y},\hat{y})$ is an approximation to  $\nabla f(\hat{y})$ given by $\nabla f(\tilde{y}) +P(\hat{y} - \tilde{y})$, with $P$ approximating the Hessian of the Lagrangian and also evaluated only at outer iterations. This strategy allows us to speed up the computations because \textit{i)} the derivatives need not be computed at each iteration \textit{ii)} efficient low-rank updates can be exploited when solving the QPs with an active set method \cite{zanelli2021}.
Notice that, unlike in \cite{Bock2007}, we do not require $P=M$.
As shown in \cite{zanelli2021},
by choosing $P$ to be the exact Hessian of the Lagrangian, we can recover local quadratic convergence, even
when $M \neq P$. In this way, we can leverage the computational
benefits
associated with choosing a positive-definite Hessian for the QP subproblems. At the same time, we recover quadratic contraction typical
of exact Hessian SQP, which generally requires a convexification strategy to handle indefinite Hessians.
\par
If convergence to a point $\tilde{z}_+ = (\tilde{y}_+, \tilde{\lambda}_+, \tilde{\mu}_+)$
is achieved, from the first-order
optimality conditions of the QPs, we obtain
\begin{equation*}\label{eq:pos}
    \begin{aligned}
         &  & a(\tilde{y},\hat{y}) + \nabla g(\tilde{y}) \tilde{\lambda}_+ + \nabla h(\tilde{y}) \tilde{\mu}_+ & = 0,                         \\
         &  & g(\tilde{y}_+)                                                                                   & = 0,                         \\
         &  & h(\tilde{y}_+)                                                                                   & \leq 0,                      \\
         &  & h_i(\tilde{y}_+) \tilde{\mu}_{+,i}                                                               & = 0, \quad i = 1,\dots, n_h,
    \end{aligned}
\end{equation*}
which can be interpreted as the first-order optimality conditions of the nonlinear program \cite{Bock2007}
\begin{equation}\label{eq:pnlp}
    \begin{aligned}
         & \underset{y}{\min} &  & f(y) + \xi(\tilde{z}, \tilde{z}_+)^{\top}y \\
         & \rm{s.t.}          &  & g(y) = 0,                                  \\
         &                    &  & h(y) \leq 0,
    \end{aligned}
\end{equation}
with
\begin{align*}
     &  & \xi(\tilde{z}, & \tilde{z}_+) \vcentcolon = \nabla f(\tilde{y}) + P(\tilde{y}_+ - \tilde{y}) - \nabla f(\tilde{y}_+) +                                                 \\
     &  &                & \left( \nabla g(\tilde{y}) - \nabla g(\tilde{y}_+)\right) \tilde{\lambda}_+ + \left( \nabla h(\tilde{y}) - \nabla h(\tilde{y}_+)\right)\tilde{\mu}_+.
\end{align*}
Since the nonlinear constraints in \eqref{eq:compactp} are preserved in the perturbed version \eqref{eq:pnlp}, the solution obtained by each outer iteration is suboptimal (due to the gradient perturbation), but feasible. The resulting strategy is described in Algorithm \ref{alg:fsqp}. The interested reader is referred to \cite{zanelli2021} for further details including globalization methods.

\begin{algorithm}
    \caption{Feasible SQP}\label{alg:fsqp}
    \begin{algorithmic}
        \Procedure{fSQP}{$\tilde{z}_0 = (\tilde{y}_0, \tilde{\lambda}_0, \tilde{\mu}_0)$, $i_{\max}$, $\epsilon_{\text{TOL}}$}
        \State $\tilde{z} \gets \tilde{z}_0$
        \For {$i=0;\,$ $i<i_{\text{max}};\,$ $i$++}
        \State compute $\nabla f(\tilde{y})$, $\nabla g(\tilde{y})$, $\nabla h(\tilde{y})$, $P$, $M$
        \State $\hat{z} \gets \tilde{z}$
        \While {$\|\Delta \hat{z}\| \geq \epsilon_{\text{TOL}}$}
        \State compute $g(\hat{y}), h(\hat{y})$ and $a(\tilde{y}, \hat{y})$
        \State solve QP (\ref{eq:compact_zo_sqp})
        \State update $\hat{y} \gets \hat{y} + \Delta y$, $\hat{\lambda} \gets \lambda_{\rm{QP}}$,
        $\hat{\mu} \gets \mu_{\rm{QP}}$
        \EndWhile
        \State $\tilde{z} \gets \hat{z}$
        \EndFor
        \State \textbf{return} $\tilde{z}$ \Comment{return feasible solution}
        \EndProcedure
    \end{algorithmic}
\end{algorithm}

\section{Inherently robust anytime MPC}\label{sec:rec_feas}
In order to be able to apply the proposed strategy to a physical system, which is inevitably affected by model mismatch, the goal is to ensure that it possesses certain robustness properties. In the following, we derive sufficient conditions for inherently robust recursive feasibility of problem \eqref{eq:nmpc} when the dynamics are affected by bounded disturbances. Notice that, since we are not using a robust-by-design approach, constraint violations may occur, but are only present when the hard-constrained version of \eqref{eq:nmpc} becomes infeasible. In that case, under the assumption that~$\mu$ is sufficiently large, known results on soft constraints \cite{Nocedal2006} guarantee that a solution with minimal constraint violation is obtained.

\subsection{Terminal Set Computation}\label{sub:terminalset}

Recursive feasibility can be achieved by constraining the terminal state to a terminal manifold $\mathcal{X}_f$ \cite{Faulwasser2009}, which is a precomputed feasible trajectory over a complete lap of the track. The terminal manifold is a positive invariant set under nominal system dynamics with a known control policy.

In this work, we use a time-varying terminal set $\mathcal{X}_f(t)$ where the terminal set at each time point is a discrete terminal state $x_f(t)$, in place of a fixed continuous terminal manifold. This greatly simplifies the design of the terminal set, and worked well in experiments (see section \ref{sec:experiments}).

The terminal states form a periodic trajectory completing one lap of the track in $T$ time steps. The terminal states at time steps $k$ and $k+1$ satisfy $x_f(k+1) = \psi(x_f(k), u_f(k), 0)$ for a known terminal control input $u_f(k)$. The terminal state trajectory $(x_f(t), u_f(t))$ is determined a priori by solving the trajectory optimization problem \eqref{eq:terminalsetoptimization}

\begin{subequations}\label{eq:terminalsetoptimization}
    \begin{align}
         & \underset{\begin{subarray}{c}
                             {x}_0, \dots, {x}_T \\
                             {u}_0, \dots, {u}_{T-1}
                         \end{subarray}}{\min} &  & \sum_{i=0}^{T-1} l({x}_i, {u}_i)                                                          \\
         & \,\,\,\quad \rm{s.t.}         &  & \hspace{-0.7ex}\left.\begin{aligned} & {x}_{i+1} - \psi({x}_i,{u}_i, 0) = 0 \\
                & (\tau_i, \delta_i) \in \mathcal{X}   \\
                & (d\tau_i, d\delta_i) \in \mathcal{U} \\
                & \pi({x}_i) \leq 0                    \\
                & i = 0,\ldots, T-1                    \\
                                                                   \end{aligned}\right\}              &                 \\
         &                               &  & {x}_T - {x}_0 = [0,0,2\pi,0,\cdots,0,l_{\text{track}}]^T,\label{eq:terminalset_periodicity}
        %2\pi\hat{\gamma} + l_\text{track}\hat{\theta}
    \end{align}
\end{subequations}
This is similar to \eqref{eq:nmpc}, but with time horizon $T$ and a periodicity constraint (\ref{eq:terminalset_periodicity}). The initial and final states are equal except that the car rotates once counterclockwise and advances forward by the length of the track~$l_{\text{track}}$.

\begin{figure}
    \centering
    \includegraphics[width=0.485\textwidth]{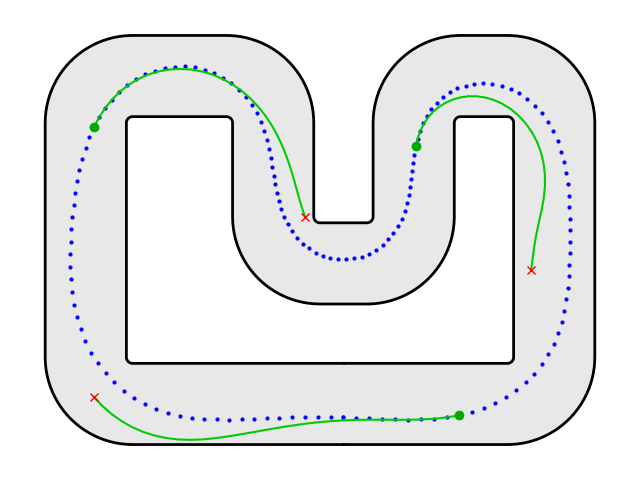}
    \caption{Terminal manifold (blue) based on a discrete feasible periodic trajectory, with example trajectories (green) satisfying the terminal constraint.}
    \label{fig:terminalset}
\end{figure}

A transitional terminal state trajectory $(\tilde{x}_f(t), \tilde{u}_f(t))$ is required for the first lap, to ensure a feasible transition from a stationary starting state to the periodic terminal manifold. We compute this by solving the similar problem (\ref{eq:terminalsetoptimization2})

\begin{subequations}\label{eq:terminalsetoptimization2}
    \begin{align}
         & \underset{\begin{subarray}{c}
                             {x}_0, \dots, {x}_{\tilde{T}} \\
                             {u}_0, \dots, {u}_{\tilde{T}-1}
                         \end{subarray}}{\min} &  & \sum_{i=0}^{\tilde{T}-1} l({x}_i, {u}_i)                                                        \\
         & \,\,\,\quad \rm{s.t.}                          &  & \hspace{-0.7ex}\left.\begin{aligned} & {x}_{i+1} - \psi({x}_i,{u}_i, 0) = 0 \\
                & (\tau_i, \delta_i) \in \mathcal{X}   \\
                & (d\tau_i, d\delta_i) \in \mathcal{U} \\
                & \pi({x}_i) \leq 0                    \\
                & i = 0,\ldots, \tilde{T}-1            \\
                                                                                    \end{aligned}\right\}   &      \\
         &                                                &  & {x}_0 = \tilde{x}_0,                                                             \\
         &                                                &  & {x}_{\tilde{T}} = x_f(0), \label{eq:terminalsetoptimization2_terminalconstraint}
    \end{align}
\end{subequations}

with a somewhat longer time horizon $\tilde{T}$ and stationary initial state $\tilde{x}_0$. The terminal constraint (\ref{eq:terminalsetoptimization2_terminalconstraint}) joins this transitional terminal trajectory to the periodic terminal trajectory computed from (\ref{eq:terminalsetoptimization}).
With this terminal set formulation, a feasible control trajectory at time step $k+1$ can be obtained from a feasible trajectory at time step $k$ by the simple shift

\begin{equation}\label{eq:shifting}
    \begin{aligned}
         & x_i(k+1) = x_{i+1}(k), \quad & i=0,\dots, N-1, \\
         & u_i(k+1) = u_{i+1}(k), \quad & i=0,\dots, N-2, \\
         & x_N(k+1) = {x_f}(k+1),       &                 \\
         & u_{N-1}(k+1) = {u_f}(k).     &                 \\
    \end{aligned}
\end{equation}

It may still be useful to explore the possibility of a continuous terminal manifold, which would allow optimizing for lap progress in real-time.

\subsection{Inherently robust recursive feasibility}
In the following, we derive the main theoretical result of the paper. We show that under mild assumptions problem~$\mathcal{P}(x,t)$ is recursively feasible with respect to the hard terminal constraint (\ref{eq:nmpc_terminal_constraint}), leading to recursive constraint satisfaction whenever the disturbance $w$ is zero and to constraint violations proportional to $w$ otherwise.

Despite a problem being feasible, a nonlinear solver may fail to converge to a solution, and otherwise is not guaranteed to do so quickly enough to be used by an MPC controller in real-time. We prove robustness to such failures by showing $M$-step open-loop recursive feasibility; under suitable assumptions, the problem remains feasible even after up to $M$ consecutive solver failures, assuming that in case of failure a new candidate is generated by shifting the current solution using (\ref{eq:shifting}).
\par
We make the following assumptions.
\begin{assumption}[Lipschitz dynamics]\label{assum:lip_dyn}
    Assume that for any $x, x', x'' \in \mathcal{X}$, any $u \in \mathcal{U}^N$ and any $w \in \mathcal{W}$, the following holds:
    \begin{equation}
        \| \phi(x'',u, 0; 1) - \phi(x',u, 0; 1)\| \leq L_{\phi,x} \cdot \| x'' - x'\|
    \end{equation}
    and
    \begin{equation}
        \| \phi(x,u, w; 1) - \phi(x,u,0; 1)\| \leq L_{\phi,w} \| w \|.
    \end{equation}
\end{assumption}

\begin{assumption}[Existence of $k$-step deadbeat controller]\label{assum:deadbeat}
    There exists a constant $\rho$ such that, for any time $t$, for any $x$ with
    $\text{dist}(x, \mathcal{X}_f(t)) \leq \rho$, there exists an input sequence $u(t) = \{u_0(t), \dots, u_{k-1}(t)\}$,
    with $u_i(t) \in \mathcal{U}, \, \forall i = 0,\dots, k-1$, such that
    $\phi(x, u(t), 0; k) \in \mathcal{X}_f(t+k)$.
\end{assumption}

\begin{figure}
    \centering
    \includegraphics[width=.485\textwidth]{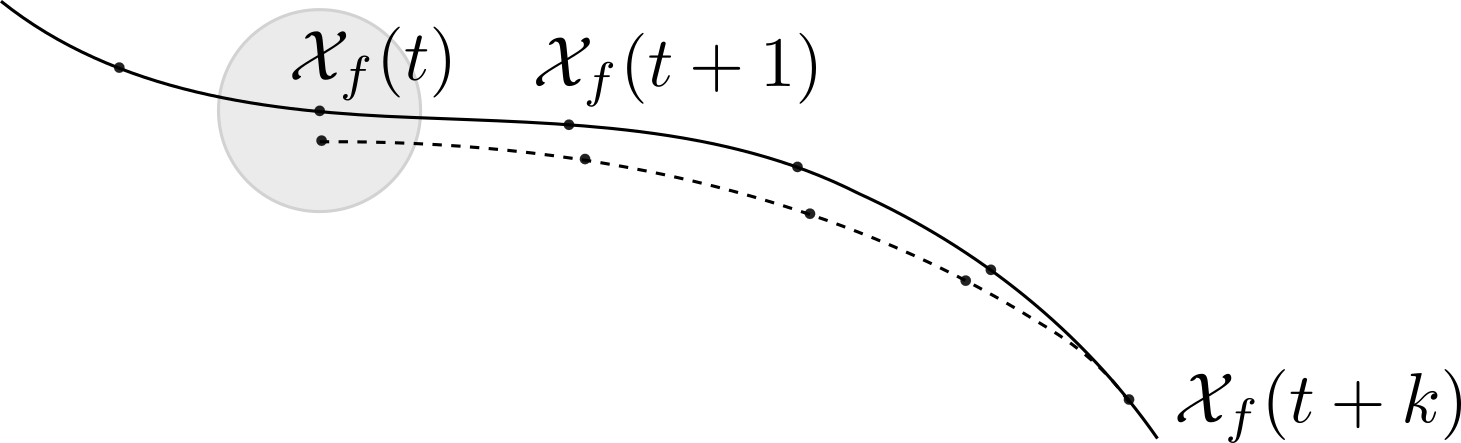}
    \caption{Illustration of time-varying terminal set and deadbeat controller.}
    \label{fig:deadbeat}
\end{figure}

Assumption \ref{assum:deadbeat} is illustrated in Figure \ref{fig:deadbeat} and, loosely speaking, requires that, if a state is sufficiently close to the time-varying terminal set it can be steered back into it in $k$ time instants.
\begin{remark}
    Designing $\mathcal{X}_f(t)$ to be a time-varying positive invariant set such that $x \in \mathcal{X}_f(t) \implies \psi(x,u_f(t), 0) \in \mathcal{X}_{f}(t+1)$, together with some regularity conditions on the reachability problem

    \begin{equation} \label{eq:reachability}
        \begin{aligned}
             &  & x             & \in \mathcal{X}_f(t),   \\
             &  & \phi(x,u,0;k) & \in \mathcal{X}_f(t+k),
        \end{aligned}
    \end{equation}
    guarantees existence of a deadbeat controller as per Assumption \ref{assum:deadbeat}. In particular, if no constraints are active and the terminal set is designed to be a single point at each time $t$, the set inclusions in (\ref{eq:reachability}) reduce to a set of nonlinear equations to be solved for the unknown input trajectory $u$.
\end{remark}

We make the following technical assumption that guarantees that for any feasible trajectory satisfying $\phi(x,u,0;N) \in \mathcal{X}_f(t)$, the nominal state $\phi(x,u,0;N+M-k)$ at time $t+M-k$ is sufficiently close to the terminal set at time $t+M-k$.
\begin{assumption}\label{assum:sigma}
    There exist a real constant $\sigma < \rho$ and an integer constant $\bar{M} \geq 1$ such that for any $x \in \mathcal{X}$, any $t$, any $M \leq \bar{M}$, and any feasible $u$ such that $\phi(x,u,0;N) \in \mathcal{X}_f(t)$ then $\text{dist}(\phi(x,u,0;N+M-k), \mathcal{X}_f(t+M-k)) < \sigma$.
\end{assumption}

\begin{remark}
    Assumption \ref{assum:sigma} can be satisfied, for example, by making sure that the discretization grid is sufficiently fine. In fact it requires that, for any trajectory ending in the terminal set $\mathcal{X}_f(t)$ associated with time $t$, the state at time $t+N+M-k$ is not too distant from $\mathcal{X}_f(t+M-k)$.
\end{remark}

We are now able to state the main theoretical result which provides robust recursive feasibility of $\mathcal{P}(x,t)$ under the assumption that the solver can fail to update the feasible candidate for at most $M$ steps.
\begin{proposition}[$M$-step open-loop recursive feasibility]\label{prop:main_prop}
    Let Assumptions \ref{assum:lip_dyn}-\ref{assum:sigma} hold. Moreover, let $x$ be such that~$\mathcal{P}(x, t)$ is feasible and let $\tilde{u}_{x}$ be a feasible solution to $\mathcal{P}(x,t)$. For any positive integer $M \geq 0$, there exists a positive constant $\delta$, such that if $w = \{w_0, w_1, \cdots\}$ satisfies $\|w_i\| \leq \delta$, $\forall i = 0, \dots, M$, then $\mathcal{P}(\phi(x,\tilde{u}_{x}, w; M), t+M)$ is feasible.
    \par
    \begin{proof}
        Since the dynamics are Lipschitz continuous, we have that, for any $x \in \mathcal{X}$ any $w \in \mathcal{W}$ and any integer $M$
        \begin{equation*}
            \| \phi(x,u, w; M) - \phi(x,u, 0; M)\| \leq \sum_{i=0}^{M-1}L_{\phi,x}^i L_{\phi,w} \delta.
        \end{equation*}
        In turn, for any $x', x''$ and any $u$, it holds that
        \begin{equation}
            \begin{aligned}
                 &  & \| \phi(x',u, 0; N-k) - & \phi(x'', u, 0; N-k)\| \leq     \\
                 &  &                         & L_{\phi,x}^{N-k} \| x'' - x'\|,
            \end{aligned}
        \end{equation}
        such that, combining the two inequalities, we obtain
        \begin{equation}
            \begin{aligned}
                 &  & \| \phi(x,u, (w^M, 0) & ; \tilde{N}) - \phi(x, u, 0; \tilde{N})\| \leq                                                                                            \\
                 &  &                       & \underbrace{L_{\phi,x}^{N-k} \sum_{i=-0}^{M-1}L_{\phi,x}^i L_{\phi,w}}_{\vcentcolon=\tilde{\sigma}} \delta = \tilde{\sigma} \cdot \delta,
            \end{aligned}
        \end{equation}
        where we have introduced $\tilde{N}\vcentcolon=N+M-k$.
        Since $\phi(x,\tilde{u}_x,0;N) \in \mathcal{X}_f(t)$, due to our assumption that $\tilde{u}_{x}$ is feasible, and $\mathcal{U}$ is bounded, together with continuity of the dynamics we have that
        \begin{equation}
            \text{dist}(\phi(x, \tilde{u}_x, 0; \tilde{N}), \mathcal{X}_f(t+M-k)) \leq \sigma
        \end{equation}
        for some bounded constant $\sigma$
        and
        \begin{equation}
            \begin{aligned}
                 &  & \text{dist}(\phi(x, \tilde{u}_x, & (w^M, 0); \tilde{N}), \mathcal{X}_f(t+M-k)) \leq \\
                 &  &                                  & \sigma + \tilde{\sigma} \delta.
            \end{aligned}
        \end{equation}
        Choosing
        \begin{equation}
            \delta \leq \frac{\rho - \sigma}{\tilde{\sigma}},
        \end{equation}
        we obtain that
        \begin{equation}
            \text{dist}(\phi(x, \tilde{u}_x, (w^M, 0);\tilde{N}), \mathcal{X}_f(t+M-k)) \leq
            \rho,
        \end{equation}
        such that we can use the deadbeat controller from Assumption \ref{assum:deadbeat} to construct a feasible candidate and $\mathcal{P}(\phi(x,\tilde{u}_x; M))$ is feasible.
    \end{proof}
\end{proposition}

\begin{figure}
    \centering
    \includegraphics[width=.45\textwidth]{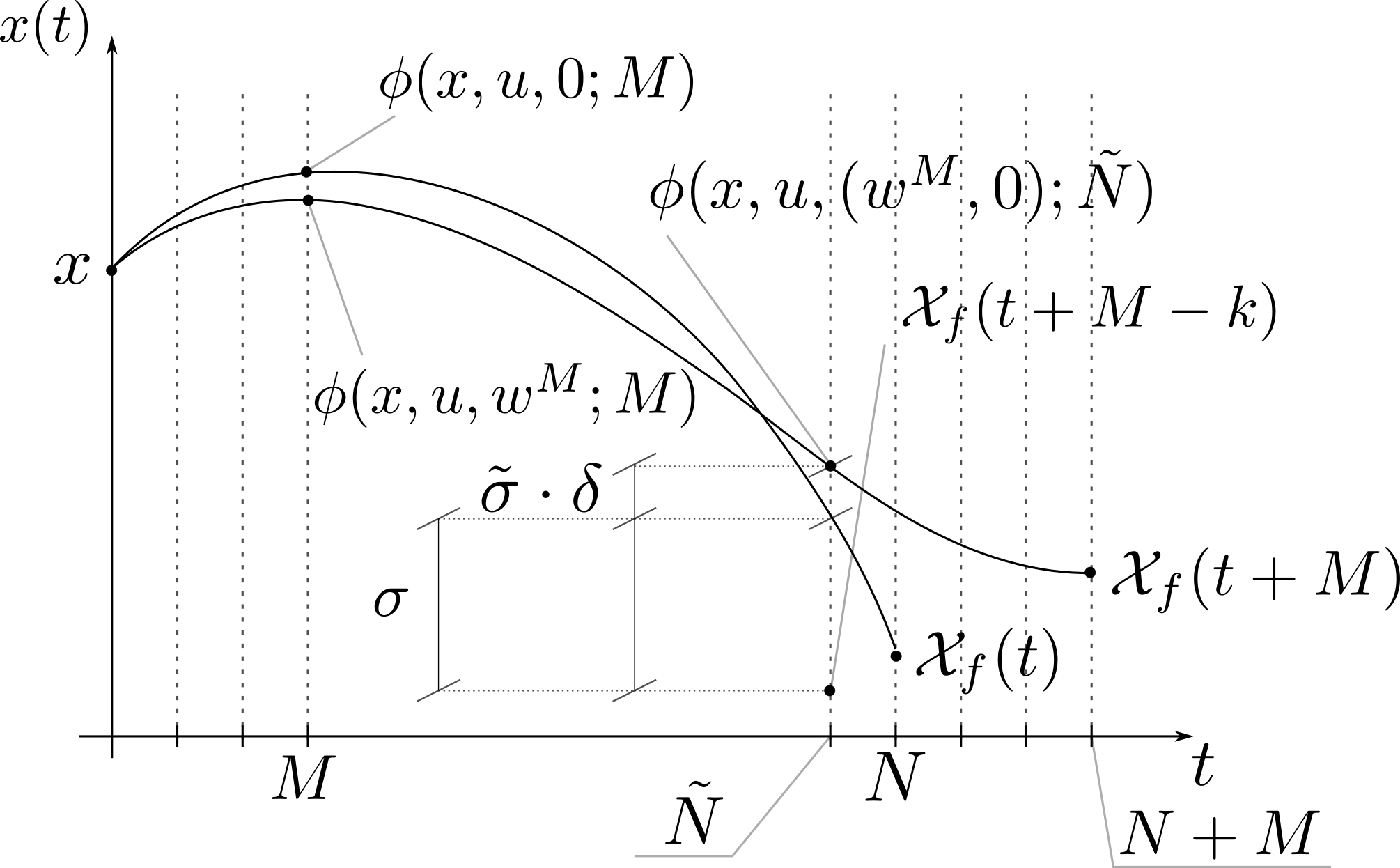}
    \caption{Illustration of the construction of the feasible candidate in the proof of Proposition \ref{prop:main_prop}. Choosing $\tilde{\sigma} \cdot \delta + \sigma \leq \rho$ we can apply the deadbeat controller from Assumption \ref{assum:deadbeat} in order to construct a feasible candidate and guarantee recursive feasibility.}
    \label{fig:proof}
\end{figure}

\section{Experiments} \label{sec:experiments}

\begin{figure}
    \centering
    \includegraphics[width=0.485\textwidth]{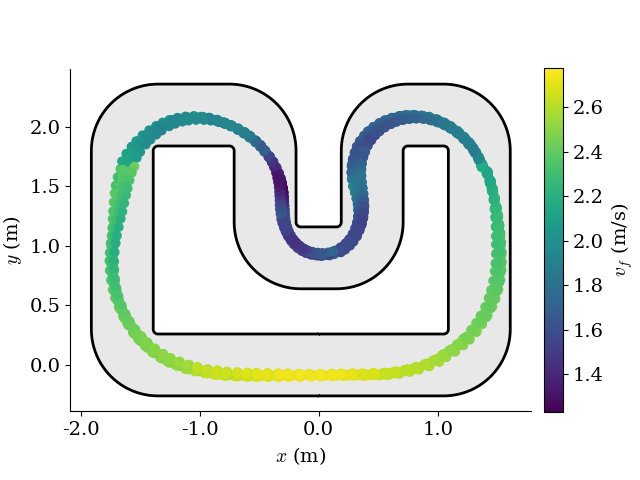}
    \caption{Race track used in simulation and hardware experiments, with closed-loop trajectories from hardware experiment HW1 colored according to forward velocity $v_f$.}
    \label{fig:closedloop}
\end{figure}

A series of experiments were performed using the feasible SQP controller to drive autonomous miniature racing vehicles around the track shown in Figure \ref{fig:closedloop}. We conducted simulation experiments to compare the closed-loop control performance of feasible SQP to other solvers, and hardware experiments to demonstrate the applicability and performance of the method in the real world. In this section, we describe the simulation and experimental setups, and then discuss the results of the experiments. All controllers run at 30 \unit{Hz} and all software runs on a Dell XPS 15 laptop using Ubuntu 20.04, Intel Core i7-9750H processor and 16 GB RAM.

\subsection{Simulation Setup}
For the simulation experiments, we use a simulation environment implemented in Python with vehicle dynamics matching the system model (section \ref{sec:systemmodel}) with additive noise. Each simulation is performed using three different controllers, each based on a different optimization solver:

\begin{itemize}
    \item FSQP, a single outer iteration of feasible SQP \\({\small\url{https://github.com/zanellia/feasible\_sqp}})
    \item RTI (real-time iteration), a single iteration of standard SQP
    \item IPOPT, an optimal interior point solver for nonconvex programming \cite{ipopt}
\end{itemize}

In each simulation, the system is run for 10 laps (1600 time steps) with each controller. Noise is applied in the form of zero-mean uniformly distributed displacement in the $p_x$ and  $p_y$ directions. The magnitude of the noise, i.e., the maximum displacement in each coordinate direction, is indicated for each simulation. The same random seed is used for all three solvers when generating noise to ensure fair comparison.

Simulation experiments are denoted SM1, SM2, SM4 and SM8, where the number indicates the maximum disturbance in each coordinate direction in \unit{\centi\metre}.

\subsection{Hardware Setup}
The hardware setup uses a Chronos miniature race car~\cite{Carron2023} (Figure \ref{fig:crs_car}) and a 46 \unit{\centi\metre} wide track with physical boundaries. A VICON motion capture system provides position and orientation measurements which are used by an extended Kalman filter to obtain complete state estimates. The software is implemented using CRS~\cite{Carron2023} in C\texttt{++}.

An important limitation of the experimental platform is that the torque output and steering effectiveness vary with the charge of the car batteries. As such, the dynamics of the vehicle are not consistent between experiments. In addition, the system model (section \ref{sec:systemmodel}) is only an approximation of the true vehicle dynamics, and some model mismatch cannot be avoided. The control problem remains recursively feasible under the inherent robustness properties described in section~\ref{sec:rec_feas}.

The car is operated in closed-loop with the feasible SQP controller. Each instance of the optimization problem is saved and later solved using RTI and IPOPT for comparison.

Hardware experiments are denoted HW1, HW2 and HW3. There is no intentional difference in configuration between hardware experiments.

\subsection{Results}\label{sub:results}

Results from both simulation and hardware experiments confirm that using the feasible SQP method, a feasible solution can be obtained with minimal additional computational effort compared to a single iteration of standard SQP which does not provide a feasible solution.

\begin{figure}
    \centering
    \includegraphics[width=.485\textwidth]{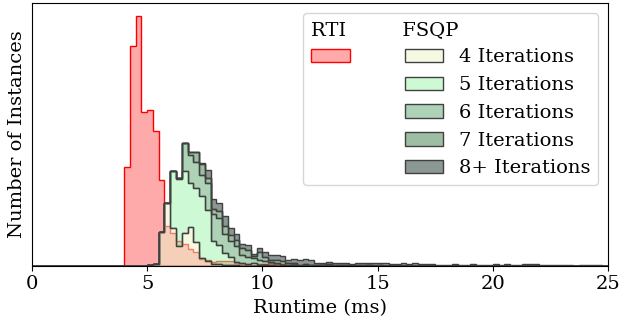}
    \caption{Histograms showing the non-normalized computation time distribution for RTI and FSQP over all optimization problems from the hardware experiments. FSQP problem instances are split by number of iterations to convergence. The weak relationship between runtime and number of iterations demonstrates that the cost of additional feasibility iterations is dominated by the cost of the first iteration (equivalent to the RTI runtime). FSQP runtime is well below 30 \unit{\milli\second}, allowing the controller to run smoothly at 30 \unit{Hz}. }
    \label{fig:iters_runtime}
\end{figure}

The feasible SQP controller ran smoothly on the hardware at the desired 30 \unit{Hz}, successfully controlling the vehicle through many consecutive laps with a competitive lap time of 5.50 seconds (165 time steps at 30 \unit{Hz}). The runtime histogram in Figure \ref{fig:iters_runtime} shows a mean computation time of 8.00 \unit{\milli\second} for our implementation of feasible SQP and that runtimes are well below 30 \unit{\milli\second} as required for a controller running at 30 \unit{Hz}. Table \ref{tbl:results} gives the average increase in runtime over standard SQP for each experiment, showing that the cost of the additional iterations to achieve feasibility is marginal.

\begin{table}[h]
    \begin{center}
        \begin{tabular}{l c c c c}
                & \multicolumn{1}{c}{FSQP}           & \multicolumn{1}{c}{FSQP/RTI} & \multicolumn{1}{c}{IPOPT/FSQP} & \multicolumn{1}{c}{FSQP/RTI} \\
                & \multicolumn{1}{c}{Convergence \%} & \multicolumn{1}{c}{Runtime}  & \multicolumn{1}{c}{Runtime}    & \multicolumn{1}{c}{Cost}     \\
                &                                    &                              &                                &                              \\
            HW1 & 97.8                               & 1.54                         & 3.57                           & 0.919                        \\
            HW2 & 98.2                               & 1.58                         & 3.58                           & 0.914                        \\
            HW3 & 97.5                               & 1.57                         & 3.76                           & 0.915                        \\
            SM1 & 98.1                               & 1.27                         & 3.88                           & 0.999998                     \\
            SM2 & 97.2                               & 1.34                         & 3.79                           & 0.9994                       \\
            SM4 & 97.7                               & 1.40                         & 3.68                           & 0.993                        \\
            SM8 & 96.0                               & 1.54                         & 3.32                           & 0.988                        \\
        \end{tabular}
    \end{center}
    \caption{Experimental results from hardware and simulation experiments. Comparison between FSQP and other solvers is averaged over cases where FSQP converged successfully. FSQP achieves consistently (though marginally) lower solution cost and produces feasible solutions with only moderately greater computation time compared to RTI and substantially faster than IPOPT. }
    \label{tbl:results}
\end{table}

\begin{figure}
    \centering
    \begin{subfigure}{.485\textwidth}
        \centering
        \includegraphics[width=\textwidth]{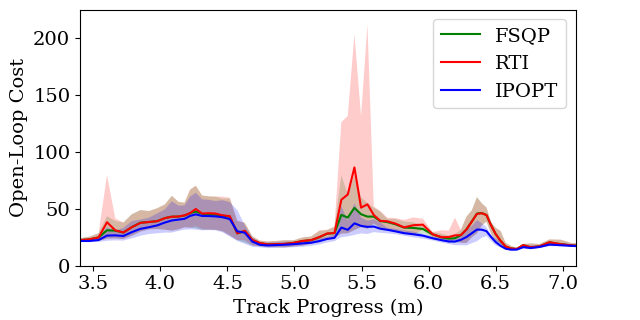}
    \end{subfigure}
    \begin{subfigure}{.485\textwidth}
        \centering
        \includegraphics[width=\textwidth]{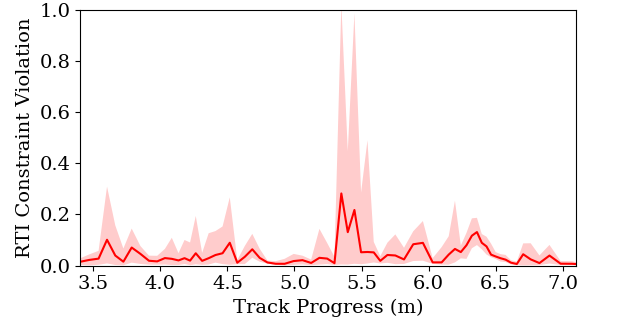}
        \caption{Hardware experiment HW1}
    \end{subfigure}
    \begin{subfigure}{.485\textwidth}
        \centering
        \includegraphics[width=\textwidth]{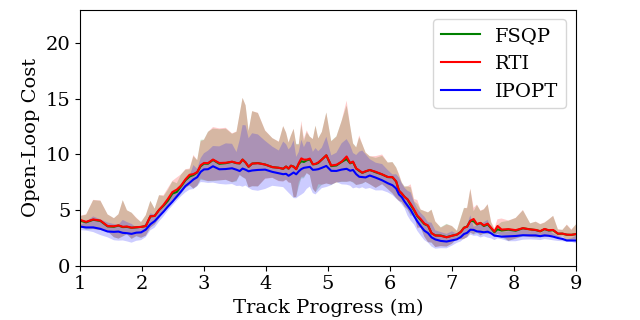}
    \end{subfigure}
    \begin{subfigure}{.485\textwidth}
        \centering
        \includegraphics[width=\textwidth]{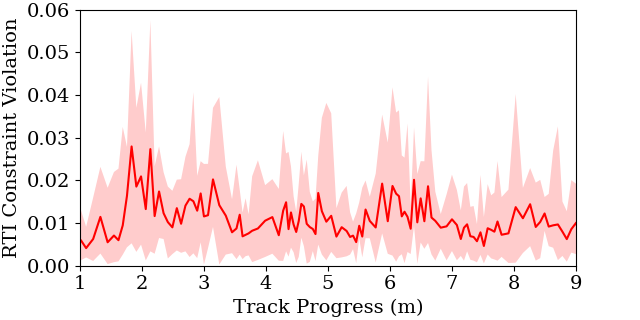}
        \caption{Simulation experiment SM4}
    \end{subfigure}
    \caption{Comparison of open-loop cost (objective function value) of solutions from FSQP, RTI and IPOPT, along with constraint violation magnitude for RTI (defined in \eqref{eq:cv}). Average values and ranges are given by plotted lines and shaded areas, respectively. FSQP solutions have comparable cost compared to RTI while eliminating constraint violation.}
    \label{fig:trajectory_quality}
\end{figure}

Figure \ref{fig:trajectory_quality} and Table \ref{tbl:results} show that the feasible SQP method achieves a moderate improvement in trajectory quality measured by objective function value compared to standard SQP, more notably in the hardware experiments. Figure \ref{fig:trajectory_quality} additionally shows the constraint violation norm for RTI, demonstrating a limitation which is solved by feasible SQP. The constraint violation is given by \eqref{eq:cv}, where $g$ and $h$ are constraint function values as in equation \eqref{eq:compactp}.

\begin{equation}\label{eq:cv}
    \| CV \|^2 \vcentcolon = \sum_{i=1}^{n_g} {g_i}^2 + \sum_{\substack{i = 1\\ h_i > 0}}^{n_h} h_i^2,
\end{equation}

It can be seen in Table \ref{tbl:results} that the implementation of feasible SQP used in this study failed to converge to a solution on a small fraction of problem instances. In these cases, the controller used the shifted solution from the previous time step, which worked well in both simulation and hardware experiments.

A video of the experiments performed can be found at:
\centerline{\url{https://youtu.be/hRQ3jOlWG6Q}}

\section{Conclusions}

In this paper, we have presented an efficient inexact MPC strategy for miniature racing using a feasible SQP algorithm which can be terminated after a single iteration while retaining recursive feasibility properties. We have demonstrated the validity of the method through experiments both in simulation and using a physical framework with autonomous miniature race cars, and have provided proof of the inherent robustness of our framework under appropriate assumptions. Comparison to state-of-the-art approximation strategies has shown that the feasible SQP algorithm provides a basis for faster MPC controllers with safety guarantees that existing methods fail to provide.

\addtolength{\textheight}{-12cm}   % This command serves to balance the column lengths

%%%%%%%%%%%%%%%%%%%%%%%%%%%%%%%%%%%%%%%%%%%%%%%%%%%%%%%%%%%%%%%%%%%%%%%%%%%%%%%%

%%%%%%%%%%%%%%%%%%%%%%%%%%%%%%%%%%%%%%%%%%%%%%%%%%%%%%%%%%%%%%%%%%%%%%%%%%%%%%%%

%%%%%%%%%%%%%%%%%%%%%%%%%%%%%%%%%%%%%%%%%%%%%%%%%%%%%%%%%%%%%%%%%%%%%%%%%%%%%%%%

\bibliographystyle{IEEEtran}
\bibliography{bibliography}

\begin{thebibliography}{10}
\providecommand{\url}[1]{#1}
\csname url@samestyle\endcsname
\providecommand{\newblock}{\relax}
\providecommand{\bibinfo}[2]{#2}
\providecommand{\BIBentrySTDinterwordspacing}{\spaceskip=0pt\relax}
\providecommand{\BIBentryALTinterwordstretchfactor}{4}
\providecommand{\BIBentryALTinterwordspacing}{\spaceskip=\fontdimen2\font plus
\BIBentryALTinterwordstretchfactor\fontdimen3\font minus
  \fontdimen4\font\relax}
\providecommand{\BIBforeignlanguage}[2]{{%
\expandafter\ifx\csname l@#1\endcsname\relax
\typeout{** WARNING: IEEEtran.bst: No hyphenation pattern has been}%
\typeout{** loaded for the language `#1'. Using the pattern for}%
\typeout{** the default language instead.}%
\else
\language=\csname l@#1\endcsname
\fi
#2}}
\providecommand{\BIBdecl}{\relax}
\BIBdecl

\bibitem{Rawlings2017}
J.~B. Rawlings, D.~Q. Mayne, and M.~M. Diehl, \emph{Model Predictive Control:
  Theory, Computation, and Design}, 2nd~ed.\hskip 1em plus 0.5em minus
  0.4em\relax Nob Hill, 2017.

\bibitem{Verschueren2014a}
R.~Verschueren, S.~D. Bruyne, M.~Zanon, J.~V. Frasch, and M.~Diehl, ``Towards
  time-optimal race car driving using nonlinear {MPC} in real-time,'' in
  \emph{CDC}, 2014, pp. 2505--2510.

\bibitem{Liniger2014}
A.~Liniger, A.~Domahidi, and M.~Morari, ``Optimization‐based autonomous
  racing of 1:43 scale {RC} cars,'' \emph{Optimal Control Applications and
  Methods}, vol.~36, pp. 628 -- 647, 2014.

\bibitem{Kabzan2020}
J.~Kabzan, L.~Hewing, A.~Liniger, and M.~N. Zeilinger, ``Learning-based model
  predictive control for autonomous racing,'' \emph{IEEE Robotics and
  Automation Letters}, vol.~4, no.~4, pp. 3363--3370, 2019.

\bibitem{Kloeser2020}
D.~Kloeser, T.~Schoels, T.~Sartor, A.~Zanelli, G.~Frison, and M.~Diehl,
  ``{NMPC} for racing using a singularity-free path-parametric model with
  obstacle avoidance,'' \emph{IFAC-PapersOnLine}, vol.~53, no.~2, pp.
  14\,324--14\,329, 2020, 21st IFAC World Congress.

\bibitem{embotech}
``embotech {AG},'' \url{www.embotech.ch}.

\bibitem{zanelli2021}
A.~Zanelli, \emph{Inexact methods for nonlinear model predictive control}, PhD
  Thesis, University of Freiburg, Germany, 2021.

\bibitem{Rajamani2006}
R.~Rajamani, \emph{Vehicle Dynamics and Control}, 2nd~ed.\hskip 1em plus 0.5em
  minus 0.4em\relax Springer US, 2012.

\bibitem{Pacejka}
H.~B. Pacejka, \emph{Tire and Vehicle Dynamics}, 3rd~ed., Elsevier Ltd., 2012.

\bibitem{Nocedal2006}
J.~Nocedal and S.~J. Wright, \emph{Numerical Optimization}, 2nd~ed., ser.
  Springer Series in Operations Research and Financial Engineering.\hskip 1em
  plus 0.5em minus 0.4em\relax Springer, 2006.

\bibitem{Allan2017}
\BIBentryALTinterwordspacing
D.~A. Allan, C.~N. Bates, M.~J. Risbeck, and J.~B. Rawlings, ``On the inherent
  robustness of optimal and suboptimal nonlinear {MPC},'' \emph{Systems \&
  Control Letters}, vol. 106, pp. 68--78, 2017. [Online]. Available:
  \url{https://www.sciencedirect.com/science/article/pii/S016769111730049X}
\BIBentrySTDinterwordspacing

\bibitem{Bock2007}
H.~G. Bock, M.~Diehl, E.~Kostina, and J.~P. Schlöder, ``Constrained optimal
  feedback control of systems governed by large differential algebraic
  equations,'' in \emph{Real-Time PDE-Constrained Optimization}.\hskip 1em plus
  0.5em minus 0.4em\relax SIAM, 2007, pp. 3--24.

\bibitem{Faulwasser2009}
T.~Faulwasser and R.~Findeisen, \emph{Nonlinear Model Predictive Path-Following
  Control}.\hskip 1em plus 0.5em minus 0.4em\relax Berlin, Heidelberg: Springer
  Berlin Heidelberg, 2009, pp. 335--343.

\bibitem{ipopt}
A.~Wächter and L.~T. Biegler, ``On the implementation of an interior-point
  filter line-search algorithm for large-scale nonlinear programming,'' in
  \emph{Mathematical Programming}, 2006, vol. 106, pp. 25--57.

\bibitem{Carron2023}
A.~Carron, S.~Bodmer, L.~Vogel, R.~Zurbrügg, D.~Helm, R.~Rickenbach,
  S.~Muntwiler, J.~Sieber, and M.~N. Zeilinger, ``Chronos and {CRS}: {D}esign
  of a miniature car-like robot and a software framework for single and
  multi-agent robotics and control,'' in \emph{2023 IEEE International
  Conference on Robotics and Automation (ICRA)}, 2023, pp. 1371--1378.

\end{thebibliography}

\end{document}